\documentclass[11pt,A4paper]{article}

\date{}
\usepackage{graphicx}
\usepackage{amsmath}
\usepackage{amsthm}
\usepackage{amssymb}
\usepackage{amstext}
\usepackage{epic}
\usepackage{eepic}
\usepackage[center]{caption2}
\usepackage{subfigure}
\usepackage{setspace}
\usepackage[noblocks]{authblk}
\allowdisplaybreaks[4]
\DeclareMathOperator{\tr}{tr}
\makeatletter
\newtheorem{theorem}[subsection]{Theorem}
\newtheorem{proposition}[subsection]{Proposition}
\newtheorem{corollary}[subsection]{Corollary}
\newtheorem{lemma}[subsection]{Lemma}
\theoremstyle{definition}
\newtheorem{definition}[subsection]{Definition}
\theoremstyle{remark}
\newtheorem{remark}[subsection]{Remark}
\makeatletter

\begin{document}
\setlength{\baselineskip}{18pt}

\title{\textbf{Some special congruences on completely regular semigroups}}
\author[1]{Li-Min Wang}
\author[2]{Ying-Ying Feng\thanks{Corresponding author. Email: rickyfungyy@fosu.edu.cn}}
\author[1]{Hong-Hua Chen}
\affil[1]{School of Mathematics, South China Normal University, \authorcr Guangzhou 510631, P. R. China}
\affil[2]{Department of Mathematics, Foshan University, \authorcr Foshan 528000, P. R. China}
\maketitle

\begin{abstract}
 This paper enriches the list of properties of the congruence sequences starting from the universal relation and successively performing the operations of lower $t$ and lower $k$. Three classes of completely regular semigroups, namely semigroups for which $\ker{\sigma}$ is a cryptogroup, semigroups for which $\ker{\nu}$ is a cryptogroup and semigroups for which $\kappa$ is over rectangular bands, are studied. $((\omega_t)_k)_t$, $((\mathcal{D}_t)_k)_t$ and $((\omega_k)_t)_k$ are found to be the least congruences on $S$ such that the quotient semigroups are semigroups for which $\ker{\sigma}$ is a cryptogroup, $\ker{\nu}$ is a cryptogroup and $\kappa$ is over rectangular bands, respectively. The results obtained present a response to three problems in Petrich and Reilly's textbook \textquoteleft\textquoteleft Completely Regular Semigroups\textquoteright\textquoteright.

 \textbf{Keywords:} completely regular semigroup, congruence, semigroup for which $\ker{\sigma}$ is a cryptogroup, semigroup for which $\ker{\nu}$ is a cryptogroup, semigroup for which $\kappa$ is over rectangular bands.

 \textbf{2010 MR Subject Classification:} 20M17
\end{abstract}

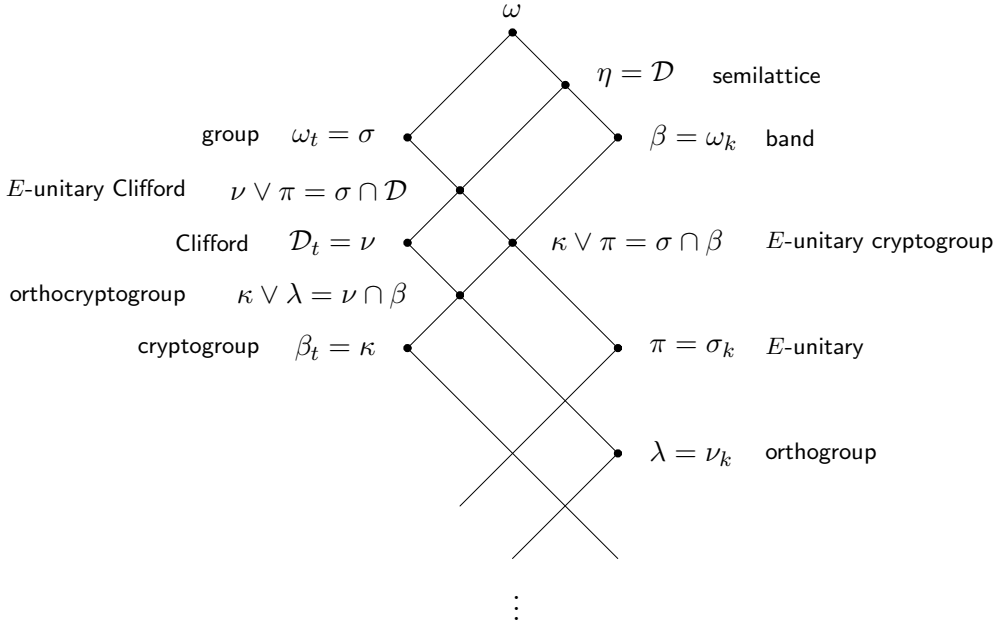
\begin{figure}[hbt]
 \renewcommand*\figurename{Figure}
 \renewcommand*\captionlabeldelim{}
 \setlength{\unitlength}{0.7cm}
 \begin{center}
 \begin{picture}(3,12)
 \drawline(2,1)(-2,5)(2,9)(0,11)(-2,9)(2,5)(-1,2)
 \drawline(0,1)(2,3)(-2,7)(1,10)
 \allinethickness{0.7mm}
 \put(-2,5){\circle*{0.06}} \put(2,5){\circle*{0.06}} \put(-2,9){\circle*{0.06}} \put(2,9){\circle*{0.06}}
 \put(0,11){\circle*{0.06}} \put(0,7){\circle*{0.06}} \put(1,10){\circle*{0.06}}
 \put(-1,8){\circle*{0.06}} \put(-2,7){\circle*{0.06}} \put(2,3){\circle*{0.06}} \put(-1,6){\circle*{0.06}}
 \put(0,11.3){\makebox(0,0)[b]{$\omega$}}
 \put(-2.6,9){\makebox(0,0)[r]{$\omega_t=\sigma$}}
 \put(2.6,9){\makebox(0,0)[l]{$\beta=\omega_k$}}
 \put(-2.6,5){\makebox(0,0)[r]{$\beta_t=\kappa$}}
 \put(2.6,5){\makebox(0,0)[l]{$\pi=\sigma_k$}}
 \put(0,0.2){\makebox(0,0)[l]{$\vdots$}}
 \put(-2,8){\makebox(0,0)[r]{$\nu \vee \pi=\sigma \cap \mathcal{D}$}}
 \put(-2,6){\makebox(0,0)[r]{$\kappa \vee \lambda=\nu \cap \beta$}}
 \put(4,7){\makebox(0,0)[r]{$\kappa \vee \pi=\sigma \cap \beta$}}
 \put(-4.8,9){\makebox(0,0)[r]{\footnotesize{\textsf{group}}}}
 \put(4.8,9){\makebox(0,0)[l]{\footnotesize{\textsf{band}}}}
 \put(-4.8,5){\makebox(0,0)[r]{\footnotesize{\textsf{cryptogroup}}}}
 \put(4.8,5){\makebox(0,0)[l]{\footnotesize{\textsf{$E$-unitary}}}}
 \put(1.6,10.2){\makebox(0,0)[l]{$\eta=\mathcal{D}$}}
 \put(-2.6,7){\makebox(0,0)[r]{$\mathcal{D}_t=\nu$}}
 \put(2.6,3){\makebox(0,0)[l]{$\lambda=\nu_k$}}
 \put(3.8,10.2){\makebox(0,0)[l]{\footnotesize{\textsf{semilattice}}}}
 \put(-5,7){\makebox(0,0)[r]{\footnotesize{\textsf{Clifford}}}}
 \put(4.8,3){\makebox(0,0)[l]{\footnotesize{\textsf{orthogroup}}}}
 \put(-6.2,8){\makebox(0,0)[r]{\footnotesize{\textsf{$E$-unitary Clifford}}}}
 \put(-6.2,6){\makebox(0,0)[r]{\footnotesize{\textsf{orthocryptogroup}}}}
 \put(4.8,7){\makebox(0,0)[l]{\footnotesize{\textsf{$E$-unitary cryptogroup}}}}
 \end{picture}
 \caption{\quad min-network of $\omega$ on completely regular semigroups} \label{original}
 \end{center}
\end{figure}

Congruences play a central role in many of the structure theorems and other important considerations in semigroup theory. The kernel-trace approach to congruences on regular semigroups is specialized to completely regular semigroups. In the attempt to gain insight into the structure of the congruence lattice $\mathcal{C}(S)$ of a completely regular semigroup $S$, a key role is played by the decomposition of $\mathcal{C}(S)$ which is induced by the complete $\cap$-congruence $\mathcal{K}$ and the complete congruence $\mathcal{T}$, induced by the kernels and the traces of congruences. For the above decomposition, we obtain two operators, lower $k$ and $t$, on $\mathcal{C}(S)$. We denote by $\rho_k$ the least congruence on $S$ having the same kernel as $\rho$, and by $\rho_t$ the least congruence on $S$ having the same trace as $\rho$. Fixing an arbitrary congruence $\rho$ on a completely regular semigroup $S$, we may apply to $\rho$ the operators of lower $t$ and $k$, thereby obtaining some new congruences. Iterating this procedure, one arrives at the min-network of congruences on $S$ based on $\rho$. This procedure becomes quite interesting if we start with the universal relation $\omega$ and in this way construct new congruences that can be characterized in an independent way. Some of these congruences are well known. Indeed, $\omega_t=\sigma$, $\omega_k=\beta$, $(\omega_t)_k=\pi$ and $(\omega_k)_t=\kappa$ are, respectively, the least group, band, $E$-unitary and cryptogroup congruences on completely regular semigroups. Also, $\mathcal{D}=\eta$, $\eta_t=\nu$ and $(\eta_t)_k=\lambda$ are, respectively, the least semilattice, Clifford and orthodox congruences on completely regular semigroups. These congruences are depicted in Figure \ref{original} together with the types of semigroups to which the quotient semigroups belong. Notice that there are corresponding results for regular semigroups (see \cite{certain}, \cite{regular}, \cite{lattice}) and inverse semigroups (see \cite{inverse}).

Recall that $\eta=\mathcal{D}$ is over groups on Clifford semigroups, and that $\beta$ is also over groups on cryptogroups. Again, notice that $\sigma$ is over rectangular bands on $E$-unitary completely regular semigroups, and $\nu$ is also over rectangular bands on orthogroups. Motivated by the symmetry observed in \cite{eomega}, we may ask: Is $\pi_t$ a congruence such that $\pi_{S/\pi_t}$ is over groups and $\lambda_t$ a congruence such that $\lambda_{S/\lambda_t}$ is over groups? Symmetrically, is $\kappa_k$ a congruence such that $\kappa_{S/\kappa_k}$ is over rectangular bands?

Coincidentally, the first question is equivalent to a problem in \cite[Problems \@Roman7.3.13]{completely}:

(\romannumeral3) For which $S$ is $\ker{\sigma}$, or $\ker{\nu}$, cryptic?\\
This problem could be reached as soon as the other two problems, especially the first one, in \cite[Problems \@Roman7.3.13]{completely} are solved:

(\romannumeral1) Characterize orthodox congruences $\rho$ for which $\ker{\rho}$ is cryptic.

(\romannumeral2) If $\mu$ is orthodox, is $\ker{\mu}$ cryptic?

Our objective here is to present a solution to these three problems. Meanwhile, we obtain properties of the congruences $\pi_t$, $\lambda_t$ and $\kappa_k$ which highlight three new classes of completely regular semigroups and lead to characterizations of these classes.

In Section \ref{pre}, we summarize the notation and terminology to be used in the paper. In Section \ref{cryptogroup}, we study semigroups for which $\ker{\sigma}$ (respectively, $\ker{\nu}$) is a cryptogroup and related congruences. A similar analysis can be found in Section \ref{over} for semigroups for which $\kappa$ is over rectangular bands and related congruences. The results obtained lead to new characterizations for orthogroups.

\section{Preliminaries}\label{pre}

Throughout this paper, unless otherwise stated, $S$ stands for a completely regular semigroup.

We will adopt the notation and terminology of Petrich -- Reilly \cite{completely} and Howie \cite{howie}, to which the reader is referred for basic information and results on completely regular semigroups. A completely regular semigroup $S$ is a union of its (maximal) subgroups. In $S$, we have the operation $a \mapsto a^{-1}$ of inversion; we will write $a^0=aa^{-1}=a^{-1}a$. For an arbitrary completely regular semigroup $S$, we denote by $E(S)$ the set of its idempotents and $V(a)$ the set of all inverses of an element $a$. The complete lattice of congruences on $S$ is denoted by $\mathcal{C}(S)$. If $S$ consists solely of idempotents, then it is a \emph{band}. A commutative band is a \emph{semilattice}. A semigroup $S$ satisfying the identity $a=axa$ is a \emph{rectangular band}. A regular semigroup is \emph{orthodox} if its idempotents form a subsemigroup. An orthodox completely regular semigroup is an \emph{orthogroup}. A semigroup $S$ is \emph{cryptic} if $\mathcal{H}$ is a congruence. A cryptic completely regular semigroup is a \emph{cryptogroup}. A semigroup isomorphic to the direct product of a rectangular band and a group is a \emph{rectangular group}.

For a regular semigroup $S$, $$a \leq b \iff a=eb=bf ~\text{for some}~ e, f \in E(S) \qquad (a, b \in S).$$ It is easy to see that the restriction of $\leq$ to $E(S)$ is given by $$e \leq f \iff e=ef=fe \qquad (e, f \in E(S)).$$

A binary relation $\theta$ on a set $X$ is denoted by $\theta_X$. For $\theta$ an equivalence relation on a semigroup $S$, we will write $\theta^*$ for the smallest congruence on $S$ containing $\theta$ and $\theta^0$ the greatest congruence contained in $\theta$. For $ Y \subseteq X$, we will write $\theta|_Y$ for \emph{restriction} of $\theta$ to $Y$. Let $\rho$ be an equivalence relation on a set $X$. If $Y \subseteq X$, then $\rho$ \emph{saturates} $Y$ if $Y$ is a union of $\rho$-classes. For a subset $K$ of a semigroup $S$, define $\pi_K$ by $$a\,\pi_K\,b \iff [\,xay \in K \iff xby \in K ~\text{for all}~ x, y \in S^1\,].$$ Then $\pi_K$ is the greatest congruence that saturates $K$.

As we've seen in the classic \cite{completely}, other than Green's relations, the relations below play an important role in the discussion of completely regular semigroups. On any completely regular semigroup $S$, relations $\mathcal{F}$, $\mathcal{Y}$ and $\Theta$ are defined by
\begin{eqnarray*}
 a\,\mathcal{F}\,b &\iff& ab^{-1} \in E(S),\\
 a\,\mathcal{Y}\,b &\iff& V(a)=V(b),\\
 a\,\Theta\,b &\iff& a^0b=ab^0.
\end{eqnarray*}

Let $\mathcal{A}$ be a class of semigroups and $\rho \in \mathcal{C}(S)$. Then $\rho$ is \emph{over} $\mathcal{A}$ if each $\rho$-class which is a subsemigroup of $S$ belongs to $\mathcal{A}$. Also $\rho$ is an $\mathcal{A}$-congruence if $S/\rho \in \mathcal{A}$. We say that $S$ is a \emph{semilattice of semigroups of type $\mathcal{P}$} if there exists a semilattice congruence on $S$ over $\mathcal{P}$. One can similarly speak of \emph{$\mathcal{C}$-bands of semigroups of type $\mathcal{P}$} if $\mathcal{C}$ is any class of bands or simply a \emph{band of semigroups of type $\mathcal{P}$} if $\mathcal{C}$ is the class of all bands. A congruence $\rho$ is \emph{idempotent separating} if $e^2=e$, $f^2=f$ and $e\,\rho\,f$ imply that $e=f$. Equivalently, $\rho$ is idempotent separating if and only if $\rho \subseteq \mathcal{H}$. On the other hand, $\rho$ is \emph{idempotent pure} if $\rho$ saturates $E(S)$; that is, $a\,\rho\,e$ with $a \in S$ and $e \in E(S)$ implies that $a \in E(S)$. Equivalently, $\rho$ is idempotent pure if and only if $\rho \subseteq \mathcal{F}$. We denote by $\mu$ and $\tau$ the greatest idempotent separating and greatest idempotent pure congruences on $S$, respectively. The equality and the universal relations on $S$ are denoted by $\varepsilon$ and $\omega$ respectively.

Some properties of completely regular semigroups are contained in the next three lemmas.

\begin{lemma}\label{green}(\cite[Lemma \@Roman2.3.5]{completely})
 Let $S$ be a regular subsemigroup of a semigroup $T$. For $\mathcal{K}=\{\mathcal{L}, \mathcal{R}, \mathcal{H}\}$, we have $\mathcal{K}_S=\mathcal{K}_T|_S$.
\end{lemma}

\begin{lemma}\label{d}(\cite[Lemma \@Roman2.3.8]{completely})
 Let $S$ be a regular subsemigroup of a completely regular semigroup $T$. Then $\mathcal{D}_S=\mathcal{D}_T|_S$.
\end{lemma}

\begin{lemma}\label{sub}(\cite[Lemma \@Roman2.3.9]{completely})
 Let $S$ be a completely regular semigroup, $e, f \in E(S)$ and $\rho \in \mathcal{C}(S)$. Then $eS$, $Sf$, $eSf$ and $e\rho$ are completely regular subsemigroups of $S$.
\end{lemma}

The following characterizations of cryptogroups, rectangular groups and orthogroups will be of use later.

\begin{lemma}\label{cryptic}(\cite[Theorem \@Roman2.8.1, Exercise \@Roman2.8.8(i)]{completely})
 The following conditions on a completely regular semigroup are equivalent.\\
 (1) $S$ is cryptic;\\
 (2) $S$ is a band of groups;\\
 (3) $S$ satisfies the identity $(ab)^0=(a^0b^0)^0$;\\
 (4) For any $a \in S$, $e \in E(S)$, $e<a^0$ implies that $ea=ae$.
\end{lemma}

\begin{lemma}\label{recg}(\cite[Corollary \@Roman3.5.3]{completely})
 The following conditions on a completely regular semigroup are equivalent.\\
 (1) $S$ is a rectangular group;\\
 (2) $E(S)$ is a rectangular band;\\
 (3) $S$ is an orthodox completely simple semigroup;\\
 (4) $\sigma$ is over rectangular bands.
\end{lemma}

\begin{lemma}\label{orthodox}(\cite[Theorem \@Roman2.5.3]{completely})
 The following conditions on a completely regular semigroup $S=(Y; S_\alpha)$ are equivalent.\\
 (1) $S$ is orthodox;\\
 (2) $S$ is a semilattice of rectangular groups;\\
 (3) For all $\alpha \in Y$, $S_\alpha$ is orthodox.
\end{lemma}
\begin{proof}
 The equivalence of (1)and (3) follows from \cite[Theorem \uppercase\expandafter{\romannumeral2}.5.3]{completely}.

 $(1) \Rightarrow (2)$. Let $S=(Y; S_\alpha)$ be an orthogroup. Then $E(S_\alpha)$ is a rectangular band by \cite[Lemma \uppercase\expandafter{\romannumeral2}.5.2]{completely}. It follows from Lemma \ref{recg} that $S_\alpha$ is a rectangular group, and $S$ is a semilattice of rectangular groups.

 $(2) \Rightarrow (1)$. If $S=(Y; S_\alpha)$ is a semilattice of rectangular groups, then $S$ is trivially a semilattice of completely simple semigroups so that $S$ is completely regular. For all $\alpha \in Y$, $S_\alpha$ is a rectangular group whence $S_\alpha$ is orthodox. That $S$ is orthodox is an immediate consequence of the equivalence of (1) and (3).
\end{proof}

We now turn to results on congruences on completely regular semigroups.

For $\rho \in \mathcal{C}(S)$, $\tr{\rho}=\rho|_{E(S)}$ is the \emph{trace} of $\rho$, and $\ker{\rho}=\{a \in S\,|\,a\,\rho\,e ~\text{for some}~ e \in E(S)\}$ is the \emph{kernel} of $\rho$. A congruence on a completely regular semigroup is determined uniquely by its trace and kernel.

\begin{lemma}\label{con}(\cite[Lemma \@Roman6.3.1]{completely})
 For any $\rho \in \mathcal{C}(S)$ and $a, b \in S$, we have $$a\,\rho\,b \iff a^0\,\tr{\rho}\,b^0 ~\text{and}~ ab^{-1} \in \ker{\rho}.$$
\end{lemma}

For any $\rho, \theta \in \mathcal{C}(S)$, the relations $\mathcal{T}$ and $\mathcal{K}$ are defined as follows, $$\rho\,\mathcal{T}\,\theta \iff \tr{\rho}=\tr{\theta}, \qquad \rho\,\mathcal{K}\,\theta \iff \ker{\rho}=\ker{\theta}.$$

The relation $\mathcal{T}$ is a complete congruence on the lattice $\mathcal{C}(S)$, while $\mathcal{K}$ is an equivalence relation on
$\mathcal{C}(S)$. The equivalence class $\rho \mathcal{T}$ [resp. $\rho \mathcal{K}$] is an interval of $\mathcal{C}(S)$ with greatest and least element to be denoted by $\rho^T$ [resp. $\rho^K$] and $\rho_t$ [resp. $\rho_k$], respectively.

\begin{lemma}\label{extreme}(\cite[Theorem \@Roman7.1.2, Theorem \@Roman7.2.2]{completely})
 For any $\rho \in \mathcal{C}(S)$, we have
 \begin{eqnarray*}
  \rho_k=(\rho \cap \mathcal{H})^*, && \rho^K=\pi_{\ker{\rho}};\\
  \rho_t=(\tr{\rho})^*=(\rho \cap \Theta)^*, && \rho^T=(\mathcal{H}(\tr{\rho})\mathcal{H})^0=(\rho \vee \mathcal{H})^0.
 \end{eqnarray*}
\end{lemma}

The following lemma will be used frequently later.

\begin{lemma}\label{theta}
 $\mathcal{H} \cap \Theta=\varepsilon$.
\end{lemma}
\begin{proof}
 Let $a, b \in S$ and $a\,(\mathcal{H} \cap \Theta)\,b$. Then $a^0=b^0$ and $a^0b=ab^0$ so that $a=aa^0=ab^0=a^0b=b^0b=b$. Hence $\mathcal{H} \cap \Theta=\varepsilon$.
\end{proof}

\section{Congruences for which $\ker{\sigma}$ of the quotient semigroup is a cryptogroup and congruences for which $\ker{\nu}$ of the quotient semigroup is a cryptogroup}\label{cryptogroup}

We first characterize orthodox congruences $\rho$ for which $\ker{\rho}$ is cryptic. Then we define semigroups for which $\ker{\sigma}$ (respectively, $\ker{\nu}$) is a cryptogroup. Characterizations for congruences for which $\ker{\sigma}$ (respectively, $\ker{\nu}$) of the quotient semigroup is a cryptogroup are considered. The least congruence $\rho$ for which $\ker{\sigma_{S/\rho}}$ is a cryptogroup and the least congruence $\theta$ for which $\ker{\nu_{S/\theta}}$ is a cryptogroup on a completely regular semigroup are proved to be $\pi_t$ and $\lambda_t$, respectively.

\begin{lemma}\label{kcg}
 The following statements concerning an orthodox congruence $\rho$ on a completely regular semigroup are equivalent.\\
 (1) $\ker{\rho}$ is cryptic;\\
 (2) $\ker{\rho} \subseteq \ker{\mu}$;\\
 (3) $\rho \cap \mathcal{H}=\rho \cap \mu$;\\
 (4) $\rho \cap \mathcal{H} \in \mathcal{C}(S)$;\\
 (5) $\rho_k \subseteq \mu$;\\
 (6) $(\rho_k)_t=\varepsilon$;\\
 (7) $\rho_k \cap \Theta=\varepsilon$;\\
 (8) $\rho_k$ is over groups.
\end{lemma}
\begin{proof}
 $(1) \Rightarrow (2)$. Let $\ker{\rho}$ be cryptic. By Lemma \ref{cryptic}, for $a \in \ker{\rho}$, $e \in E(\ker{\rho})=E(S)$, $e \leq a^0$ implies $ea=ae$, and therefore by \cite[Corollary \@Roman6.2.9]{completely}, $a \in \ker{\mu}$. Hence $\ker{\rho} \subseteq \ker{\mu}$.

 $(2) \Rightarrow (3)$. Let $a, b \in S$ with $a\,(\rho \cap \mathcal{H})\,b$. Then $a^0=b^0$ and $ab^{-1} \in \ker{\rho} \subseteq \ker{\mu}$ whence $a\,\mu\,b$ by Lemma \ref{con}. Hence $\rho \cap \mathcal{H} \subseteq \rho \cap \mu$. The opposite conclusion is trivial whence equality prevails.

 $(3) \Rightarrow (4)$. Obvious.

 $(4) \Rightarrow (5)$. From Lemma \ref{extreme} and the assumption we see that $\rho_k=(\rho \cap \mathcal{H})^*=\rho \cap \mathcal{H} \subseteq \mathcal{H}$. Consequently $\rho_k$ is idempotent separating so that $\rho_k \subseteq \mu$.

 $(5) \Rightarrow (6)$. Since $\rho_k \subseteq \mu$, it follows that $(\rho_k)_t \subseteq \mu_t=\varepsilon$, i.e. that $(\rho_k)_t=\varepsilon$.

 $(6) \Rightarrow (7)$. Indeed $\tr{\rho_k}=\tr{(\rho_k)_t}=\varepsilon$ so that $\rho_k$ is idempotent separating whence $\rho_k \subseteq \mu$. Thus $\rho_k=(\rho \cap \mathcal{H})^* \subseteq \mu \subseteq \mathcal{H}$ and $\rho_k \subseteq \rho \cap \mathcal{H}$. The reverse containment is obvious. Therefore we have $\rho_k=\rho \cap \mathcal{H}$. By Lemma \ref{theta} we can deduce that $\rho_k \cap \Theta=\rho \cap \mathcal{H} \cap \Theta=\rho \cap \varepsilon=\varepsilon$.

 $(7) \Rightarrow (8)$. Suppose that $e, f \in E(S)$ with $e\,\rho_k\,f$. Since any two idempotents are $\Theta$-related, we have $e\,(\rho_k \cap \Theta)\,f$ which yields $e=f$. Therefore $\rho_k \subseteq \mathcal{H}$. Since $e\rho_k$ is a completely regular subsemigroup of $S$ by Lemma \ref{sub}, it follows that $e\rho_k$ is a subgroup of $H_e$.

 $(8) \Rightarrow (1)$. Notice first that $\ker{\rho}$ is a subsemigroup of $S$ because $\rho$ is an orthodox congruence. By \cite[Lemma \@Roman6.2.6]{completely}, the hypothesis implies that $\rho_k \subseteq \mu \subseteq \mathcal{H}$. It follows that $\rho_k=\rho \cap \mathcal{H}$ from the proof of $(5) \Rightarrow (6)$ and $(6) \Rightarrow (7)$, and so, $\rho_k|_{\ker{\rho}} \subseteq \mathcal{H}|_{\ker{\rho}}$. For any $a, b \in \ker{\rho}$, if $a\,\mathcal{H}\,b$, then $a\,\rho\,a^0=b^0\,\rho\,b$ and $a\,(\rho \cap \mathcal{H})\,b$. Hence $a\,\rho_k\,b$ and $\mathcal{H}|_{\ker{\rho}} \subseteq \rho_k|_{\ker{\rho}}$. Therefore $\mathcal{H}|_{\ker{\rho}}=\rho_k|_{\ker{\rho}}$. Note that $\rho_k|_{\ker{\rho}}$ is a congruence on $\ker{\rho}$. Hence $\mathcal{H}_{\ker{\rho}}=\mathcal{H}|_{\ker{\rho}}$ is a congruence on $\ker{\rho}$ which implies that $\ker{\rho}$ is cryptic.
\end{proof}

\begin{remark}
 It is worth noticing that the assumption that $\rho$ is orthodox is only used in part (1) --- otherwise $\ker{\rho}$ might not be a subsemigroup --- and proving the implication $(8) \Rightarrow (1)$.
\end{remark}

\begin{corollary}\label{mu}
 If $\mu$ is orthodox, then $\ker{\mu}$ is cryptic.
\end{corollary}
\begin{proof}
 If $\mu$ is orthodox, then from the equivalence of (1) and (2) in Lemma \ref{kcg} we infer that $\ker{\mu}$ is cryptic.
\end{proof}

\begin{remark}
 Lemma \ref{kcg} and Corollary \ref{mu} are exactly the responses to Problems \@Roman7.3.13 (\romannumeral1) and (\romannumeral2) in \cite{completely}.
\end{remark}

We next proceed to answer the question in \cite[Problems \@Roman7.3.13(\romannumeral3)]{completely}.

\begin{theorem}\label{sigmacg}
 The following conditions on a completely regular semigroup $S$ are equivalent.\\
 (1) $\ker{\sigma}$ is cryptic;\\
 (2) $\ker{\sigma} \subseteq \ker{\mu}$;\\
 (3) $\sigma \cap \mathcal{H}=\sigma \cap \mu$;\\
 (4) $\sigma \cap \mathcal{H} \in \mathcal{C}(S)$;\\
 (5) $\pi \subseteq \mu$;\\
 (6) $\pi_t=\varepsilon$;\\
 (7) there exists an idempotent separating $E$-unitary congruence on $S$;\\
 (8) $\pi \cap \Theta=\varepsilon$;\\
 (9) $\pi$ is over groups.
\end{theorem}
\begin{proof}
 Notice that we have used no assumption of orthodoxy within the proof of $(2) \Rightarrow (3) \Rightarrow (4) \Rightarrow (5) \Rightarrow (6) \Rightarrow (7) \Rightarrow (8)$ in Lemma \ref{kcg}. So we now have $(2) \Rightarrow (3) \Rightarrow (4) \Rightarrow (5) \Rightarrow (6) \Rightarrow (8) \Rightarrow (9)$. It suffices to show $(1) \Rightarrow (2)$, $(9) \Rightarrow (1)$ and the equivalence of (6) and (7).

 $(1) \Rightarrow (2)$. Notice that the assumption that $\rho$ is orthodox is a guarantee for $\ker{\rho}$ being a subsemigroup in Lemma \ref{kcg}. Here we have $\ker{\sigma}$ being a cryptogroup. And the argument goes along the same lines as in $(1) \Rightarrow (2)$ of Lemma \ref{kcg}.

 $(9) \Rightarrow (1)$. First it follows from \cite[Theorem \@Roman7.3.6]{completely} that $\ker{\sigma}=\ker{\pi}$ is a subsemigroup of $S$. A similar argument to that for $(8) \Rightarrow (1)$ in Lemma \ref{kcg} establishes that $\mathcal{H}_{\ker{\sigma}}$ is a congruence on $\ker{\sigma}$ and so the proof is complete.

 $(6) \Rightarrow (7)$. It is clear from the hypothesis that $\tr{\pi}=\tr{\pi_t}=\varepsilon$, and that the $E$-unitary congruence $\pi$ is idempotent separating.

 $(7) \Rightarrow (6)$. Let $\rho$ be an idempotent separating $E$-unitary congruence. Then $\pi \subseteq \rho \subseteq \mu$ and $\pi_t \subseteq \mu_t=\varepsilon$ which gives $\pi_t=\varepsilon$.
\end{proof}

We now set out to set up characterizations for congruences for which $\ker{\sigma}$ of the quotient semigroup is a cryptogroup. We shall need an auxiliary result.

For a congruence $\rho$ on a completely regular semigroup $S$, notice that there exists at least one cryptogroup congruence on $S$ containing $\rho$, namely the universal relation $\omega$, and that the intersection of cryptogroup congruences is still a cryptogroup congruence, we have a least cryptogroup congruence containing $\rho$, to be denoted by $\kappa_\rho$. Similarly, we have a least $E$-unitary and a least orthodox congruence on $S$ containing $\rho$, to be denoted by $\pi_\rho$ and $\lambda_\rho$, respectively.

\begin{proposition}\label{rho}
 Let $S$ be a completely regular semigroup and $\rho \in \mathcal{C}(S)$. Then\\
 (1) $\kappa_\rho=\rho \vee \kappa$, and $\kappa_{S/\rho}=\kappa_\rho/\rho=(\rho \vee \kappa)/\rho$;\\
 (2) $\pi_\rho=((\rho \vee \sigma)\cap \rho\mathcal{H}\rho)^*$, and $\pi_{S/\rho}=\pi_\rho/\rho$;\\
 (3) $\lambda_\rho=((\rho \vee \nu)\cap \rho\mathcal{H}\rho)^*$, and $\lambda_{S/\rho}=\lambda_\rho/\rho$.
\end{proposition}
\begin{proof}
 (1) First let $a, b \in S$. Recalling Lemma \ref{cryptic}, we have $((a\kappa)(b\kappa))^0=((a\kappa)^0(b\kappa)^0)^0$ which implies that $(ab)^0\kappa=((ab)\kappa)^0=((a^0\kappa)(b^0\kappa))^0=(a^0b^0)^0\kappa$. Further, $\kappa \subseteq \rho \vee \kappa$ implies that $(ab)^0(\rho \vee \kappa)=(a^0b^0)^0(\rho \vee \kappa)$ and thus also that $((a(\rho \vee \kappa))(b(\rho \vee \kappa)))^0=((a(\rho \vee \kappa))^0(b(\rho \vee \kappa))^0)^0$. It now follows from Lemma \ref{cryptic} that $\rho \vee \kappa$ is a cryptogroup congruence.

 Next let $\theta$ be a cryptogroup congruence containing $\rho$. Then $\kappa \subseteq \theta$ and $\rho \subseteq \theta$. Therefore $\rho \vee \kappa \subseteq \theta$ proving the minimality of $\rho \vee \kappa$.

 Now $(S/\rho)/(\kappa_\rho/\rho) \simeq S/\kappa_\rho$. Thus $\kappa_\rho/\rho$ is a cryptogroup congruence on $S/\rho$. If $\theta/\rho$ is a cryptogroup congruence on $S/\rho$ with $\rho \subseteq \theta$, then $S/\theta \simeq (S/\rho)/(\theta/\rho)$ which implies that $\theta$ is a cryptogroup congruence on $S$. Hence $\kappa_\rho \subseteq \theta$ and $\kappa_\rho/\rho \subseteq \theta/\rho$. Consequently $\kappa_\rho/\rho$ is the least cryptogroup congruence on $S/\rho$ whence $\kappa_{S/\rho}=\kappa_\rho/\rho$. Notice that similar arguments give that $\pi_{S/\rho}=\pi_\rho/\rho$ and $\lambda_{S/\rho}=\lambda_\rho/\rho$.

 (2) The remark above shows that for any $a, b \in S$, $$a\,\pi_\rho\,b \iff (a\rho)\,\pi_{S/\rho}\,(b\rho).$$ By Lemma \ref{extreme} we have an expression for $\pi_{S/\rho}$, hence it suffices to prove $$a\,((\rho \vee \sigma) \cap \rho\mathcal{H}\rho)^*\,b \iff (a\rho)\,(\sigma_{S/\rho} \cap \mathcal{H}_{S/\rho})^*\,(b\rho),$$ and for this, it is enough to show $$a\,((\rho \vee \sigma) \cap \rho\mathcal{H}\rho)\,b \iff (a\rho)\,(\sigma_{S/\rho} \cap \mathcal{H}_{S/\rho})\,(b\rho).$$ It is further enough to show the equivalences
 \begin{eqnarray*}
  a\,(\rho \vee \sigma)\,b &\iff& (a\rho)\,\sigma_{S/\rho}\,(b\rho),\\
  a\,(\rho\mathcal{H}\rho)\,b &\iff& (a\rho)\,\mathcal{H}_{S/\rho}\,(b\rho).
 \end{eqnarray*}
 The first equivalence follows immediately from the fact that $\rho \vee \sigma$ is the least group congruence containing $\rho$. By \cite[Theorem \@Roman6.5.1]{completely}, the second equivalence holds. We consequently have $\pi_\rho=((\rho \vee \sigma)\cap \rho\mathcal{H}\rho)^*$.

 (3) The proof is closely similar to that for (2) and is omitted.
\end{proof}

We are now ready for characterizations for congruences for which $\ker{\sigma}$ of the quotient semigroup is a cryptogroup.

\begin{proposition}\label{sigmacgc}
 The following statements concerning $\rho \in \mathcal{C}(S)$ are equivalent.\\
 (1) $\ker{\sigma_{S/\rho}}$ is cryptic;\\
 (2) $\pi_\rho \subseteq \rho^T$;\\
 (3) $\tr{\pi_\rho}=\tr{\rho}$.
\end{proposition}
\begin{proof}
 $(1) \Rightarrow (2)$. Since $\ker{\sigma_{S/\rho}}$ is a cryptogroup, we have by Theorem \ref{sigmacg} that $\pi_{S/\rho} \subseteq \mu_{S/\rho}$. It is known by Proposition \ref{rho} that $\pi_{S/\rho}=\pi_\rho/\rho$ and by \cite[Proposition \@Roman7.2.5]{completely} that $\mu_{S/\rho}=\rho^T/\rho$. Hence $\pi_\rho/\rho \subseteq \rho^T/\rho$ so that $\pi_\rho \subseteq \rho^T$.

 $(2) \Rightarrow (3)$. $\pi_\rho \subseteq \rho^T$ gives $\tr{\pi_\rho} \subseteq \tr{\rho^T}=\tr{\rho}$. But $\rho \subseteq \pi_\rho$, and then $\tr{\rho} \subseteq \tr{\pi_\rho}$. Consequently, $\tr{\pi_\rho}=\tr{\rho}$.

 $(3) \Rightarrow (1)$. The hypothesis and $\rho \subseteq \pi_\rho$ give that the $E$-unitary congruence $\pi_\rho/\rho$ is idempotent separating on $S/\rho$. By Theorem \ref{sigmacg}, we get that $\ker{\sigma_{S/\rho}}$ is a cryptogroup.
\end{proof}

\begin{corollary}\label{pit}
 $\pi_t=((\omega_t)_k)_t$ is the least congruence $\rho$ for which $\ker{\sigma_{S/\rho}}$ is a cryptogroup on a completely regular semigroup.
\end{corollary}
\begin{proof}
 Since $\pi_{\pi_t}=\pi \subseteq \pi^T=(\pi_t)^T$, we have by Proposition \ref{sigmacgc} that $\ker{\sigma_{S/\pi_t}}$ is a cryptogroup. If $\rho$ is a congruence for which $\ker{\sigma_{S/\rho}}$ is a cryptogroup, then $\pi \subseteq \pi_\rho \subseteq \rho^T$ and $\pi_t \subseteq (\rho^T)_t=\rho_t \subseteq \rho$. This proves the assertion.
\end{proof}

We conclude this section with observations concerning semigroups for which $\ker{\nu}$ is a cryptogroup.

\begin{theorem}\label{nucg}
 The following conditions on a completely regular semigroup $S$ are equivalent.\\
 (1) $\ker{\nu}$ is cryptic;\\
 (2) $\ker{\nu} \subseteq \ker{\mu}$;\\
 (3) $\nu \cap \mathcal{H}=\nu \cap \mu$;\\
 (4) $\nu \cap \mathcal{H} \in \mathcal{C}(S)$;\\
 (5) $\lambda \subseteq \mu$;\\
 (6) $\lambda_t=\varepsilon$;\\
 (7) there exists an idempotent separating orthodox congruence on $S$;\\
 (8) $\lambda \cap \Theta=\varepsilon$;\\
 (9) $\lambda$ is over groups.
\end{theorem}
\begin{proof}
 First notice that $\nu$ is orthodox from $\nu\,\mathcal{K}\,\lambda$ and \cite[Exercises \@Roman7.3.11(\romannumeral1)]{completely}. So the equivalence of (1) -- (6), (8) and (9) is an immediate consequence of Lemma \ref{kcg}. It remains to show the equivalence of (6) and (7).

 $(6) \Rightarrow (7)$. It is clear from the hypothesis that $\tr{\lambda}=\tr{\lambda_t}=\varepsilon$, and that the orthodox congruence $\lambda$ is idempotent separating.

 $(7) \Rightarrow (6)$. Let $\rho$ be an idempotent separating orthodox congruence. Then $\lambda \subseteq \rho \subseteq \mu$ and $\lambda_t \subseteq \mu_t=\varepsilon$ which establishes $\lambda_t=\varepsilon$.
\end{proof}

The next result characterizes congruences for which $\ker{\nu}$ of the quotient semigroup is a cryptogroup in several ways.

\begin{proposition}\label{nucgc}
 The following statements concerning $\rho \in \mathcal{C}(S)$ are equivalent.\\
 (1) $\ker{\nu_{S/\rho}}$ is cryptic;\\
 (2) $\lambda_\rho \subseteq \rho^T$;\\
 (3) $\tr{\lambda_\rho}=\tr{\rho}$.
\end{proposition}
\begin{proof}
 The proof is closely similar to that for Proposition \ref{sigmacgc} and is omitted.
\end{proof}

\begin{corollary}
 $\lambda_t=((\mathcal{D}_t)_k)_t$ is the least congruence $\rho$ for which $\ker{\nu_{S/\rho}}$ is a cryptogroup on a completely regular semigroup.
\end{corollary}
\begin{proof}
 The argument here goes along the same lines as in Corollary \ref{pit}.
\end{proof}

\begin{remark}
 Theorem \ref{sigmacg} and Theorem \ref{nucg} are the responses to Problem \@Roman7.3.13(\romannumeral3) in \cite{completely}.
\end{remark}

\section{Congruences for which $\kappa$ of the quotient semigroup is over rectangular bands}\label{over}

After defining semigroups for which $\kappa$ is over rectangular bands, we provide some equivalent conditions in terms of congruences. We then characterize congruences for which $\kappa$ of the quotient semigroup is over rectangular bands on a completely regular semigroup $S$ and prove that the least such congruence on $S$ is $\kappa_k$. These results lead naturally to congruence characterizations for orthogroups.

\begin{definition}\label{defko}
 A completely regular semigroup $S$ is called a \emph{semigroup for which $\kappa$ is over rectangular bands} if $e\kappa$ is a rectangular band for each $e \in E(S)$.
\end{definition}

The next proposition illustrates this class of completely regular semigroups.

\begin{proposition}\label{brecg}
 Let $S$ be a semigroup for which $\kappa$ is over rectangular bands. Then $S$ is a band of rectangular groups.
\end{proposition}
\begin{proof}
 By \cite[Lemma \@Roman6.1.9]{completely}, $\kappa \cap \mathcal{H}=\varepsilon$, and therefore, $\kappa$ is idempotent pure, i.e., $\ker{\kappa}=E(S)$. Since $\kappa \subseteq \beta$, we must have $e\kappa \subseteq e\beta$ for all $e \in E(S)$ whence $e\beta$ is a completely regular subsemigroup of $S$ by Lemma \ref{sub}. From $\kappa\,\mathcal{T}\,\beta$ and the hypothesis that $e\kappa$ is a rectangular band, we now infer that $e\beta$ is a rectangular group by Lemma \ref{recg}. Thus $S$ is a band of rectangular groups.
\end{proof}

In fact, if $S$ is a band of rectangular groups and $\rho$ a band congruence on $S$ such that $e\rho$ is a rectangular group for all $e \in E(S)$, then $\beta \subseteq \rho$ and $e\beta$ is a completely regular subsemigroup of $e\rho$ by Lemma \ref{sub}. It follows by Lemma \ref{recg} that $e\beta$ is a rectangular group. For convenience, we write $\underset{\alpha \in B}{\bigcup} S_\alpha$ or $(B; S_\alpha)$ for $S=\underset{e \in E(S)}{\bigcup} e\beta$, where $B$ is a band and $S_\alpha$ is a rectangular group.

The next lemma provides useful information concerning bands of rectangular groups.

\begin{lemma}\label{ftheta}
 Let $S=(B; S_\alpha)$ be a band of rectangular groups. Then $\beta \cap \mathcal{F}=\beta \cap \Theta$.
\end{lemma}
\begin{proof}
 Let $a, b \in S$ with $a\,(\beta \cap \mathcal{F})\,b$. Then $ab^{-1} \in E(S)$ and $a, b \in S_\alpha$ for some $\alpha \in B$. It is clear that $ab^{-1} \in S_\alpha$ whence $b^0=b^0ab^{-1}b^0$. It follows that $b=b^0ab^0$ so that $a^0b=a^0b^0ab^0=a^0b^0a^0ab^0=a^0ab^0=ab^0$. Hence $a\,\Theta\,b$ and $\beta \cap \mathcal{F} \subseteq \beta \cap \Theta$. On the other hand, if $a\,(\beta \cap \Theta)\,b$, then $ab^0=a^0b$ which implies $ab^{-1}=a^0bb^{-1}=a^0b^0 \in E(S_\alpha)$. Thus $a\,\mathcal{F}\,b$ and $\beta \cap \Theta \subseteq \beta \cap \mathcal{F}$. Therefore we have $\beta \cap \mathcal{F}=\beta \cap \Theta$.
\end{proof}

It is worth remarking that orthogroups are semigroups for which $\kappa$ is over rectangular bands. First, if $S$ is an orthogroup, then $\lambda=\varepsilon$ and $\kappa \subseteq \tau$, which yields that $e\kappa$ is a band for $e \in E(S)$. Next recall that an orthodox completely regular semigroup is a semilattice of rectangular groups. For $f, g \in e\kappa$, $f, g \in S_\alpha$ since $\kappa \subseteq \beta \subseteq \mathcal{D}$. But $S_\alpha$ is a rectangular group and $f=fgf$ and hence $e\kappa$ is a rectangular band.

The class of semigroups for which $\kappa$ is over rectangular bands can be characterized in many ways, of which we now give a sample.

\begin{theorem}\label{korecb}
 The following conditions on a completely regular semigroup $S$ are equivalent.\\
 (1) $\kappa$ is over rectangular bands;\\
 (2) $\beta \cap \mathcal{F}=\beta \cap \Theta=\beta \cap \tau$;\\
 (3) $\beta \cap \Theta \in \mathcal{C}(S)$;\\
 (4) $\kappa \cap \mathcal{H}=\varepsilon$;\\
 (5) $\kappa \subseteq \tau$;\\
 (6) $\kappa_k=\varepsilon$;\\
 (7) there exists an idempotent pure cryptogroup congruence on $S$;\\
 (8) $\kappa \subseteq \mathcal{F}$;\\
 (9) $S$ is a band of rectangular groups with $\tr{\beta} \subseteq \tr{\tau}$;\\
 (10) $\kappa \subseteq \mathcal{Y}$.
\end{theorem}
\begin{proof}
 $(1) \Rightarrow (2)$. Let $\kappa$ be over rectangular groups. By Proposition \ref{brecg}, $S$ is a band of rectangular groups and so we can apply Proposition \ref{ftheta} to conclude that $\beta \cap \mathcal{F}=\beta \cap \Theta$. Since $\kappa$ is over rectangular bands, we have $\ker{\kappa}=E(S)$. Moreover, $\ker{\kappa}=\ker{\beta}$ yields $\beta \subseteq \tau$ and $\beta \cap \mathcal{F} \subseteq \beta=\beta \cap \tau$. Conversely, $\tau \subseteq \mathcal{F}$ gives $\beta \cap \tau \subseteq \beta \cap \mathcal{F}$. And the required equality holds.

 $(2) \Rightarrow (3)$. Obvious.

 $(3) \Rightarrow (4)$. It follows from $\kappa=\beta_t=(\beta \cap \Theta)^*=\beta \cap \Theta$ and Lemma \ref{theta} that $\kappa \cap \mathcal{H}=\beta \cap \Theta \cap \mathcal{H}=\beta \cap \varepsilon=\varepsilon$.

 $(4) \Rightarrow (5)$. It follows directly from \cite[Lemma \@Roman6.1.8]{completely}.

 $(5) \Rightarrow (6)$. If $\kappa \subseteq \tau$, then $\kappa_k \subseteq \tau_k=\varepsilon$ whence $\kappa_k=\varepsilon$.

 $(6) \Rightarrow (7)$. If $\kappa_k=\varepsilon$, then $\ker{\kappa}=\ker{\kappa_k}=E(S)$ which implies that the cryptogroup congruence $\kappa$ is idempotent pure.

 $(7) \Rightarrow (8)$. Let $\rho$ be an idempotent pure cryptogroup congruence. Then $\kappa \subseteq \rho \subseteq \tau$. If $a, b \in S$ and $a\,\kappa\,b$, then $ab^{-1}\,\kappa\,b^0$. From the foregoing we have $ab^{-1} \in E(S)$, that is, $a\,\mathcal{F}\,b$, and $\kappa \subseteq \mathcal{F}$.

 $(8) \Rightarrow (1)$. Again notice that $e\kappa$ is a completely regular subsemigroup of $S$ by Lemma \ref{sub}. Let $a \in S$, $e \in E(S)$ and $a\,\kappa\,e$. Then $a^0\,\kappa\,e\,\kappa\,a$. By the hypothesis this implies that $a\,\mathcal{F}\,a^0$ whence $a=aa^0=a(a^0)^{-1} \in E(S)$. Therefore $e\kappa$ is a band. Further, $\kappa \subseteq \beta \subseteq \mathcal{D}$ implies that $e\kappa$ is completely simple. We conclude that $e\kappa$ is a rectangular band.

 $(1) \Rightarrow (9)$. Let $\kappa$ be over rectangular bands. It follows from Proposition \ref{brecg} that $S$ is a band of rectangular groups and from the previous proof of $(1) \Rightarrow (5)$ that $\kappa \subseteq \tau$ whence $\tr{\beta}=\tr{\kappa} \subseteq \tr{\tau}$.

 $(9) \Rightarrow (2)$. Let $a, b \in S$ with $a\,(\beta \cap \mathcal{F})\,b$. Then $a^0\,\tr{\beta}\,b^0$. But $\tr{\beta} \subseteq \tr{\tau}$, and hence $a^0\,\tr{\tau}\,b^0$. It follows from Lemma \ref{con} that $a\,\tau\,b$ since $ab^{-1} \in E(S)$. Thus $\beta \cap \mathcal{F} \subseteq \beta \cap \tau$. Conversely, if $a\,(\beta \cap \tau)\,b$, then $ab^{-1}\,\tau\,b^0$ which implies that $ab^{-1} \in E(S)$, i.e., $a\,\mathcal{F}\,b$. Therefore $\beta \cap \tau \subseteq \beta \cap \mathcal{F}$ whence $\beta \cap \mathcal{F}=\beta \cap \tau$. Again, by Proposition \ref{ftheta}, we have $\beta \cap \mathcal{F}=\beta \cap \Theta$.

 $(1) \Leftrightarrow (10)$. It follows immediately from \cite[Lemma \@Roman6.1.9]{completely}.
\end{proof}

We now turn to characterizations of congruences for which $\kappa$ of the quotient semigroup is over rectangular bands on a completely regular semigroup.

\begin{proposition}\label{korecbc}
 The following statements concerning $\rho \in \mathcal{C}(S)$ are equivalent.\\
 (1) $\kappa_{S/\rho}$ is over rectangular bands;\\
 (2) $\rho^K$ is a cryptogroup congruence;\\
 (3) $\ker{(\rho \vee \kappa)}=\ker{\rho}$.
\end{proposition}
\begin{proof}
 $(1) \Rightarrow (2)$. Since $\kappa_{S/\rho}$ is over rectangular bands, we have by Theorem \ref{korecb} that $\kappa_{S/\rho} \subseteq \tau_{S/\rho}$. It is known by Proposition \ref{rho} that $(\rho \vee \kappa)/\rho \subseteq \rho^K/\rho$ so that $\rho \vee \kappa \subseteq \rho^K$. Therefore $\kappa \subseteq \rho^K$ and $\rho^K$ is a cryptogroup congruence since cryptogroups are closed under homomorphic images.

 $(2) \Rightarrow (3)$. It follows from $\kappa \subseteq \rho^K$ that $\rho \vee \kappa \subseteq \rho^K$ whence $\ker{\rho} \subseteq \ker{(\rho \vee \kappa)} \subseteq \ker{\rho^K}=\ker{\rho}$ and $\ker{(\rho \vee \kappa)}=\ker{\rho}$.

 $(3) \Rightarrow (1)$. Since $\ker{(\rho \vee \kappa)}=\ker{\rho}$, $(\rho \vee \kappa)/\rho$ is an idempotent pure congruence on $S/\rho$, and this together with Proposition \ref{rho} gives that $\kappa_{S/\rho}$ is over rectangular bands by Theorem \ref{korecb}.
\end{proof}

\begin{corollary}
 In any completely regular semigroup $S$, $\kappa_k=((\omega_k)_t)_k$ is the least congruence $\rho$ for which $\kappa_{S/\rho}$ is over rectangular bands.
\end{corollary}
\begin{proof}
 The argument here goes along the same lines as in Corollary \ref{nucgc}.
\end{proof}

\begin{proposition}
 For any completely regular semigroup $S$, $\kappa \cap \pi=\pi_t \vee \kappa_k$ is the least congruence $\rho$ on $S$ for which $\ker{\sigma_{S/\rho}}$ is a cryptogroup and $\kappa_{S/\rho}$ is over rectangular bands.
\end{proposition}
\begin{proof}
 By \cite[Proposition \@Roman7.2.10]{completely}, we have $\kappa \cap \pi=(\sigma \cap \beta)_t \cap (\sigma \cap \beta)_k=(\sigma \cap \beta)_{tk} \vee (\sigma \cap \beta)_{kt}=\kappa_k \vee \pi_t$. Note that $\pi_{\kappa \cap \pi}=\pi=\sigma_k \subseteq \omega_k=\beta$. Hence $\tr{(\kappa \cap \pi)}=\tr{\kappa} \cap \tr{\pi}=\tr{\beta} \cap \tr{\pi}=\tr{\pi}$, and so $\ker{\sigma_{S/(\kappa \cap \pi)}}$ is a cryptogroup, by Proposition \ref{sigmacgc}. Also, $\ker{(\kappa \cap \pi)}=\ker{\kappa} \cap \ker{\pi}=\ker{\kappa} \cap \ker{\sigma}=\ker{\kappa}=\ker{((\pi \cap \kappa)\vee \kappa)}$, which implies that $\kappa_{S/(\kappa \cap \pi)}$ is over rectangular bands. Let $\rho$ be a congruence such that $\ker{\sigma_{S/\rho}}$ is a cryptogroup and that $\kappa_{S\rho}$ is over rectangular bands. By Proposition \ref{korecbc} and Proposition \ref{sigmacgc}, we obtain $\kappa \subseteq \rho^K$ and $\pi \subseteq \pi_{\rho} \subseteq \rho^T$, which yields $\kappa \cap \pi \subseteq \rho^K \cap \rho^T=\rho$. This completes the proof of the proposition.
\end{proof}

\begin{figure}[hbt]
 \renewcommand*\figurename{Figure}
 \renewcommand*\captionlabeldelim{}
 \setlength{\unitlength}{0.7cm}
 \begin{center}
 \begin{picture}(6,16)
 \drawline(-1.5,1.5)(2,5)(-2,9)(2,13)(0,15)(-2,13)(2,9)(-2,5)(1,2)
 \drawline(-0.5,1.5)(-2,3)(2,7)(-2,11)(1,14)
 \allinethickness{0.7mm}
 \put(-2,5){\circle*{0.06}} \put(2,5){\circle*{0.06}}
 \put(-2,9){\circle*{0.06}} \put(2,9){\circle*{0.06}}
 \put(-2,13){\circle*{0.06}} \put(2,13){\circle*{0.06}}
 \put(0,15){\circle*{0.06}} \put(0,11){\circle*{0.06}} \put(0,7){\circle*{0.06}} \put(0,3){\circle*{0.06}}
 \put(1,14){\circle*{0.06}} \put(-1,12){\circle*{0.06}} \put(-2,11){\circle*{0.06}} \put(2,7){\circle*{0.06}} \put(-2,3){\circle*{0.06}} \put(-1,10){\circle*{0.06}}
 \put(0,15.3){\makebox(0,0)[b]{$\omega$}}
 \put(-2.6,13){\makebox(0,0)[r]{$\omega_t=\sigma$}}
 \put(2.6,13){\makebox(0,0)[l]{$\beta=\omega_k$}}
 \put(-2.6,9){\makebox(0,0)[r]{$\beta_t=\kappa$}}
 \put(2.6,9){\makebox(0,0)[l]{$\pi=\sigma_k$}}
 \put(-2.6,5){\makebox(0,0)[r]{$\pi_t$}}
 \put(2.6,5){\makebox(0,0)[l]{$\kappa_k$}}
 \put(0,0.7){\makebox(0,0)[l]{$\vdots$}}
 \put(-4.8,13){\makebox(0,0)[r]{\footnotesize{\textsf{group}}}}
 \put(4.8,13){\makebox(0,0)[l]{\footnotesize{\textsf{band}}}}
 \put(-4.8,9){\makebox(0,0)[r]{\footnotesize{\textsf{cryptogroup}}}}
 \put(4.8,9){\makebox(0,0)[l]{\footnotesize{\textsf{$E$-unitary}}}}
 \put(-3.7,5){\makebox(0,0)[r]{\footnotesize{\textsf{$\ker{\sigma}$ is cryptic}}}}
 \put(3.7,5){\makebox(0,0)[l]{\footnotesize{\textsf{$\kappa$ is over rectangular bands}}}}
 \put(1.6,14.2){\makebox(0,0)[l]{$\eta=\mathcal{D}$}}
 \put(-2.6,11){\makebox(0,0)[r]{$\mathcal{D}_t=\nu$}}
 \put(2.6,7){\makebox(0,0)[l]{$\lambda=\nu_k$}}
 \put(-2.6,3){\makebox(0,0)[r]{$\lambda_t$}}
 \put(3.8,14.2){\makebox(0,0)[l]{\footnotesize{\textsf{semilattice}}}}
 \put(-5,11){\makebox(0,0)[r]{\footnotesize{\textsf{Clifford}}}}
 \put(4.8,7){\makebox(0,0)[l]{\footnotesize{\textsf{orthogroup}}}}
 \put(-3.8,3){\makebox(0,0)[r]{\footnotesize{\textsf{$\ker{\nu}$ is cryptic}}}}
 \end{picture}
 \caption{\quad The least congruences for which $\ker{\sigma_{S/\rho}}$ is a cryptogroup, $\ker{\nu_{S/\rho}}$ is a cryptogroup and $\kappa_{S/\rho}$ is over rectangular bands} \label{new}
 \end{center}
\end{figure}
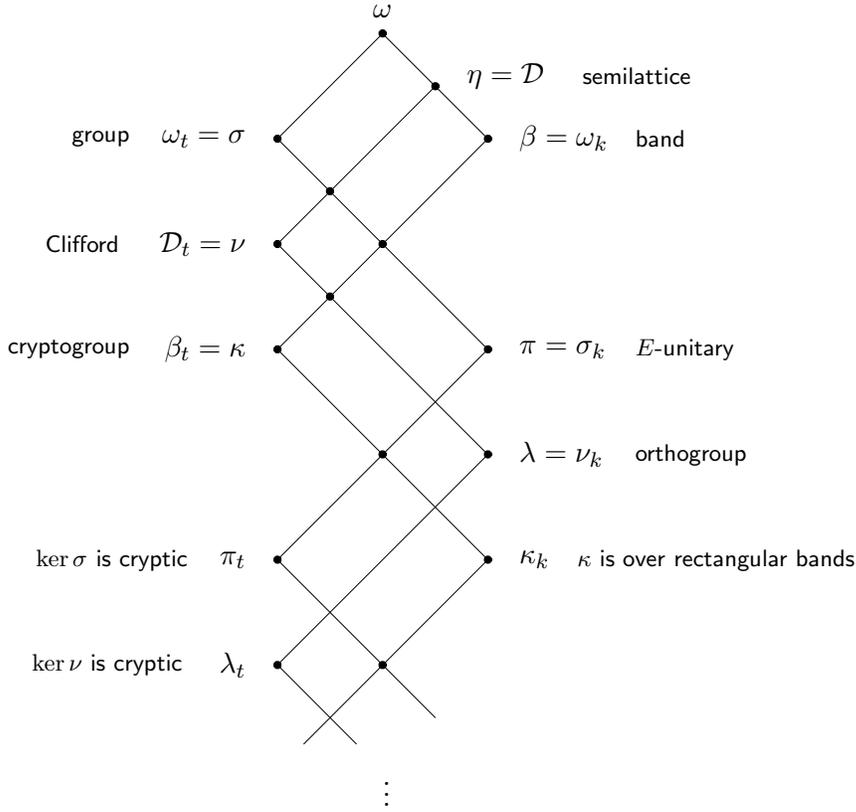

The congruences $\pi_t$, $\lambda_t$ and $\kappa_k$ are depicted in Figure \ref{new} together with some well-known congruences, enlarging the list of properties of the congruence network in Figure \ref{original}.

The next proposition explores an interesting property of semigroups for which $\kappa$ is over rectangular bands.

\begin{proposition}
 Let $S=\underset{\alpha \in B}{\bigcup}S_\alpha$ be a semigroup for which $\kappa$ is over rectangular bands and $\sigma_\alpha$ the least group congruence on $S_\alpha$ for any $\alpha \in B$. Then $\kappa=\underset{\alpha \in B}{\bigcup}\sigma_\alpha$.
\end{proposition}
\begin{proof}
 First since $\kappa\,\mathcal{T}\,\beta$ and since $\kappa \subseteq \beta$, we have $\tr{\kappa|_{S_\alpha}}=\omega|_{E(S_\alpha)}$. Thus $\kappa|_{S_\alpha}$ is a group congruence on $S_\alpha$ and $\sigma_\alpha \subseteq \kappa|_{S_\alpha}$. Conversely, suppose that $a, b \in S_\alpha$ and $a\,\kappa|_{S_\alpha}\,b$. Then $ab^{-1}\,\kappa\,b^0$. By the hypothesis that $\kappa$ is over rectangular bands, $ab^{-1} \in E(S)$ and therefore $b=b^0b\,\sigma_\alpha\,ab^{-1}b=ab^0\,\sigma_\alpha\,aa^0=a$. Consequently $\kappa|_{S_\alpha} \subseteq \sigma_\alpha$ so that $\kappa|_{S_\alpha}=\sigma_\alpha$.

 Next, it is clear that $\kappa|_{S_\alpha} \subseteq \kappa$ for all $\alpha \in B$. Then $\underset{\alpha \in B}{\bigcup}\kappa|_{S_\alpha} \subseteq \kappa$. Conversely, if $a, b \in S$ and $a\,\kappa\,b$, then $a\,\beta\,b$ so that $a$, $b$ lies in some $S_\alpha$. Thus $a\,\kappa|_{S_\alpha}\,b$. We conclude that $\kappa=\underset{\alpha \in B}{\bigcup}\kappa|_{S_\alpha}=\underset{\alpha \in B}{\bigcup}\sigma_\alpha$.
\end{proof}

Recall that orthogroups are examples of semigroups for which $\kappa$ is over rectangular bands. The remainder of this section is devoted to new characterizations for orthogroups.

\begin{lemma}\label{ot}
 Let $S$ be an orthogroup. Then $\rho_t=(\rho \cap \mathcal{F})^*$.
\end{lemma}
\begin{proof}
 Since $\rho \cap \mathcal{F} \subseteq \rho$, $(\rho \cap \mathcal{F})^* \subseteq \rho$ and thereby $\tr{(\rho \cap \mathcal{F})^*} \subseteq \tr{\rho}$. On the other hand, for any $e, f \in E(S)$, $e\,\mathcal{F}\,f$ since $S$ is orthodox. Therefore $\tr{\rho} \subseteq \tr{(\rho \cap \mathcal{F})} \subseteq \tr{(\rho \cap \mathcal{F})^*}$ so that $\tr{(\rho \cap \mathcal{F})^*}=\tr{\rho}$.

 Let $\theta \in \mathcal{C}(S)$ and $\rho\,\mathcal{T}\,\theta$. If $a, b \in S$ and $a\,(\rho \cap \mathcal{F})\,b$, then $a\,\rho\,b$ and $ab^{-1} \in E(S)$. It follows that $a^0\,\rho\,b^0$ while $a^0\,\theta\,b^0$. But $ab^{-1} \in E(S) \subseteq \ker{\theta}$ and hence $a\,\theta\,b$, by Lemma \ref{con}. Therefore $\rho \cap \mathcal{F} \subseteq \theta$ whence $(\rho \cap \mathcal{F})^* \subseteq \theta$. We conclude that $\rho_t=(\rho \cap \mathcal{F})^*$.
\end{proof}

Compare the following result with Theorem \ref{korecb} and Proposition \ref{brecg}.

\begin{theorem}\label{orthogroup}
 The following conditions on a completely regular semigroup are equivalent.\\
 (1) $S$ is orthodox;\\
 (2) $S$ is a semilattice of rectangular groups;\\
 (3) $\mathcal{D} \cap \mathcal{F}=\mathcal{D} \cap \tau$;\\
 (4) $\mathcal{D} \cap \mathcal{F} \in \mathcal{C}(S)$;\\
 (5) $\mathcal{D} \cap \Theta=\mathcal{D} \cap \tau$;\\
 (6) $\mathcal{D} \cap \Theta \in \mathcal{C}(S)$;\\
 (7) $\nu \subseteq \tau$;\\
 (8) $\nu \cap \mathcal{H}=\varepsilon$;\\
 (9) $\nu_k=\varepsilon$;\\
 (10) there exists an idempotent pure Clifford congruence;\\
 (11) $\nu \subseteq \mathcal{F}$;\\
 (12) $\nu$ is over rectangular bands;\\
 (13) $\nu \subseteq \mathcal{Y}$.
\end{theorem}
\begin{proof}
 $(1) \Leftrightarrow (2)$. See Lemma \ref{orthodox}.

 $(1) \Rightarrow (3)$. Assume that $a, b \in S$ and $a\,(\mathcal{D} \cap \mathcal{F})\,b$. Then $a\,\mathcal{D}\,b$. From the equivalence of (1) and (2), we know that $E(D_{a^0})=E(D_{b^0})$ is a rectangular band. Therefore by \cite[Exercises \@Roman2.5.9(\romannumeral1)]{completely}, if $x, y \in S^1$ and $xa^0y \in E(S)$, then $xb^0(a^0y) \in E(S)$. Again, by \cite[Exercises \@Roman2.5.9(\romannumeral1)]{completely}, $(xb^0a^0)b^0y \in E(S)$, i.e., $xb^0y \in E(S)$. Similarly, if $xb^0y \in E(S)$, then $xa^0y \in E(S)$. Therefore $a^0\,\tau\,b^0$. Furthermore, $a\,\mathcal{F}\,b$ yields $ab^{-1} \in E(S)$. By Lemma \ref{con}, we can deduce that $a\,\tau\,b$ whence $\mathcal{D} \cap \mathcal{F} \subseteq \mathcal{D} \cap \tau$. Conversely if $a\,(\mathcal{D} \cap \tau)\,b$, then $ab^{-1}\,\tau\,b^0$ which implies that $ab^{-1} \in E(S)$, i.e., $a\,\mathcal{F}\,b$. Thus $\mathcal{D} \cap \tau \subseteq \mathcal{D} \cap \mathcal{F}$ and the required equality holds.

 $(3) \Rightarrow (4)$, $(5) \Rightarrow (6)$. Obvious.

 $(4) \Rightarrow (7)$. Suppose that $\mathcal{D} \cap \mathcal{F}$ is a congruence. Then $\mathcal{D} \cap \mathcal{F} \subseteq \mathcal{F}$ gives $\mathcal{D} \cap \mathcal{F} \subseteq \tau$, by \cite[Lemma \@Roman6.1.8]{completely}. Hence $\mathcal{D} \cap \mathcal{F} \subseteq \mathcal{D} \cap \tau$. It is clear that $\tau \subseteq \mathcal{F}$. Thus $\mathcal{D} \cap \tau \subseteq \mathcal{D} \cap \mathcal{F}$ and the equality $\mathcal{D} \cap \mathcal{F}=\mathcal{D} \cap \tau$ takes place. From the foregoing and Lemma \ref{ot} we infer that $\nu=\mathcal{D}_t=(\mathcal{D} \cap \mathcal{F})^*=\mathcal{D} \cap \mathcal{F}=\mathcal{D} \cap \tau \subseteq \tau$.

 $(2) \Rightarrow (5)$. As a first step we show that $\mathcal{D} \cap \mathcal{F}=\mathcal{D} \cap \Theta$. Let $a, b \in S$ and $a\,(\mathcal{D} \cap \mathcal{F})\,b$. Then $a$, $b$ lies in some rectangular group $S_\alpha$ and $ab^{-1} \in E(S_\alpha)$. Thus $b^0=b^0ab^{-1}b^0=b^0ab^{-1}$, while $b=b^0ab^0$. It follows that $a^0b=a^0b^0ab^0=a^0b^0a^0ab^0=a^0ab^0=ab^0$, i.e., $a\,\Theta\,b$, and so $\mathcal{D} \cap \mathcal{F} \subseteq \mathcal{D} \cap \Theta$. Conversely, if $a\,(\mathcal{D} \cap \Theta)\,b$, then $ab^0=a^0b$ and $ab^{-1}=a^0b^0 \in E(S)$. Thus $a\,\mathcal{F}\,b$ and $\mathcal{D} \cap \Theta \subseteq \mathcal{D} \cap \mathcal{F}$. Therefore $\mathcal{D} \cap \mathcal{F}=\mathcal{D} \cap \Theta$.

 From the fact that $(1) \Rightarrow (3)$ and the equivalence of (1) and (2), we have $\mathcal{D} \cap \Theta=\mathcal{D} \cap \mathcal{F}=\mathcal{D} \cap \tau$.

 $(6) \Rightarrow (7)$. It follows from $\nu=\mathcal{D}_t=(\mathcal{D} \cap \Theta)^*=\mathcal{D} \cap \Theta$ and Lemma \ref{theta} that $\nu \cap \mathcal{H}=\mathcal{D} \cap \Theta \cap \mathcal{H}=\mathcal{D} \cap \varepsilon=\varepsilon$.

 $(7) \Rightarrow (8)$, $(11) \Leftrightarrow (12) \Leftrightarrow (13)$. It follows immediately from Lemma \@Roman6.1.8 and Lemma \@Roman6.1.9 in \cite{completely}.

 $(8) \Rightarrow (9)$. From $\nu \cap \mathcal{H}=\varepsilon$ we infer that $\nu_k=(\nu \cap \mathcal{H})^*=\varepsilon^*=\varepsilon$.

 $(9) \Rightarrow (10)$. Indeed $\ker{\nu}=\ker{\nu_k}=E(S)$ so that the Clifford congruence $\nu$ is idempotent pure.

 $(10) \Rightarrow (11)$. Let $\rho$ be an idempotent pure Clifford congruence. Then $\nu \subseteq \rho \subset \tau$, which by \cite[Lemma \@Roman6.1.8]{completely} gives $\nu \subseteq \mathcal{F}$.

 $(12) \Rightarrow (2)$. Suppose that $\nu$ is over rectangular bands. Since $\mathcal{D}\,\mathcal{T}\,\nu$, this implies, in view of Lemma \ref{recg}, that $e\mathcal{D}$ is a rectangular group for any $e \in E(S)$. Therefore, $S$ is a semilattice of rectangular groups.
\end{proof}

\textbf{Acknowledgements} This work is supported by the National Natural Science Foundation of China (Grant No. 11871150). The authors would like to thank the referee for the careful reading and valuable suggestions. The authors would also like to thank Professor Marianne Johnson for providing the efficient communications.


\begin{thebibliography}{99}
 \bibitem{howie} Howie, J. M. (1995). \textit{Fundamentals of Semigroup Theory, London Mathematical Society Monographs}. New Series, Vol.12. New York: The Clarendon Press.
 \bibitem{certain} Howie, J. M., Lallement, G. (1966). Certain fundamental congruences on a regular semigroup. \textit{Proc. Glasgow Math. Assoc.} 7: 145--159.
 \bibitem{regular} Pastijn, F., Petrich, M. (1986). Congruences on regular semigroups. \textit{Trans. Amer. Math. Soc.} 295: 607--633 .
 \bibitem{lattice} Pastijn, F., Petrich, M. (1988). The congruence lattice of a regular semigroup. \textit{J. Pure and Appl. Algebra} 53: 93--123.
 \bibitem{inverse} Petrich, M. (1984). \textit{Inverse Semigroups}. New York: Wiley.
 \bibitem{completely} Petrich, M., Reilly, N. R. (1999). \textit{Completely Regular Semigroups}. New York: Wiley.
 \bibitem{eomega} Wang, L. M., Feng, Y. Y. (2011). $E\omega$-Clifford congruences and $E\omega$-$E$-reflexive congruences on an inverse semigroup. \textit{Semigroup Forum} 82: 354--366. DOI: 10.1007/s00233-011-9293-y.
\end{thebibliography}
\end{document}